**Konstantine Zelator**

# Heron isosceles Triangles with integral external radii $\rho_\alpha, \rho_\beta, \rho_\gamma$

Department of Mathematics,
College of Arts and Sciences,
 Mail stop 942
The University of Toledo,
Toledo, OH 43606-3390
e-mails: konstantine-zelator@yahoo.com
         konstantine.zelator@utoledo.edu

**April 26, 2007**

# 1 Introduction

For each triangle on a Euclidean plane, there exist three characteristic exterior circles. Each of these circles is tangent to one of the triangle's sides; and also tangential to the two straight lines containing the other two sides, but not the sides themselves. Each of the three circle centers is the point of intersection between two of the triangles' external angle bisectors as well as one internal angle bisector (figure 1). In a recently published paper by Amy Bell ([1]), a number of properties of the radii, $\rho_\alpha, \rho_\beta, \rho_\gamma$ respectively say, for such circles, are explored. There is also the work of Hansen ([2]) in which results pertaining to these radii in the case of right angled triangles are obtained.

Throughout this paper we will denote by $\alpha, \beta, \gamma$ the three side lengths of a triangle $ABC$.

The purpose of this work is two-fold. The first objective is to present an alternative simple derivation of the formulas for the radii $\rho_\alpha, \rho_\beta, \rho_\gamma$; in terms of triangle $ABC$'s side lengths $\alpha, \beta, \gamma$. We achieve this in Section 3. To do so, we employ the Law of Cosines and two simple trigonometric identities. For a geometric proof (rather than one which involves trigonometric functions), the reader may refer to [6]. For related material see [7], [8], [9]. In Section 4, we state the well known parametric formulas which describe the entire family of Pythagorean triangles or triples. In Section 5, we combine the formulas from Sections 3 and 4. This we do, in order to give formulas for the radii $\rho_\alpha, \rho_\beta, \rho_\gamma$, for triangles which are Pythagorean. In Section 6, we examine the family of all isosceles triangles with integer side lengths and integral area. This is the main or primary goal of this work; to parametrically describe (in terms of three integer-valued parameters) the set of all Heron (see below) isosceles triangles which have integral radii $\rho_\alpha, \rho_\beta, \rho_\gamma$. We interject here a note on terminology.

Note that, typically, a triangle whose side lengths and area are integers is referred to as a **Heron triangle.** At times we will be using this terminology. Also note that an isosceles **Heron triangle** must obviously be non equilateral; if $\alpha = \beta = \gamma$ is the integer side length of an equilateral triangle its area $E$ is given by $E = \dfrac{\alpha^2 \cdot \sqrt{3}}{4}$, an irrational number.

References for Heron triangles can be found in [10] and [11].

In the section 6, we use an auxiliary proposition, Proposition 1, in order to establish Proposition 2, which provides us with a parametric description of all Heron isosceles triangles. Furthermore, in Theorem 1, we parametrically describe a subset $S$ of all Heron isosceles triangles, namely, those Heron isosceles triangles for which $\rho_\alpha, \rho_\beta, \rho_\gamma$ are also integers. In the last section, Section 7, we present numerical examples of such triangles by means of two tables.

## 2 The Law of Cosines and two trigonometric identities

With reference to Figure 1, $\theta$ denotes the interior angle $A$ of triangle $ABC$ and according to the Law of Cosines, $\alpha^2 = \beta^2 + \gamma^2 - 2\beta\gamma\cos\theta$. Thus,

$$\cos\theta = \frac{\beta^2 + \gamma^2 - \alpha^2}{2\beta\gamma}$$

(1)

Moreover,

$$\cos\theta = 2\cos^2(\frac{\theta}{2}) - 1 = 1 - 2\sin^2(\frac{\theta}{2});$$

which gives

$$\tan\left(\frac{\theta}{2}\right) = \sqrt{\frac{1 - \cos\theta}{1 + \cos\theta}} \qquad (2)$$

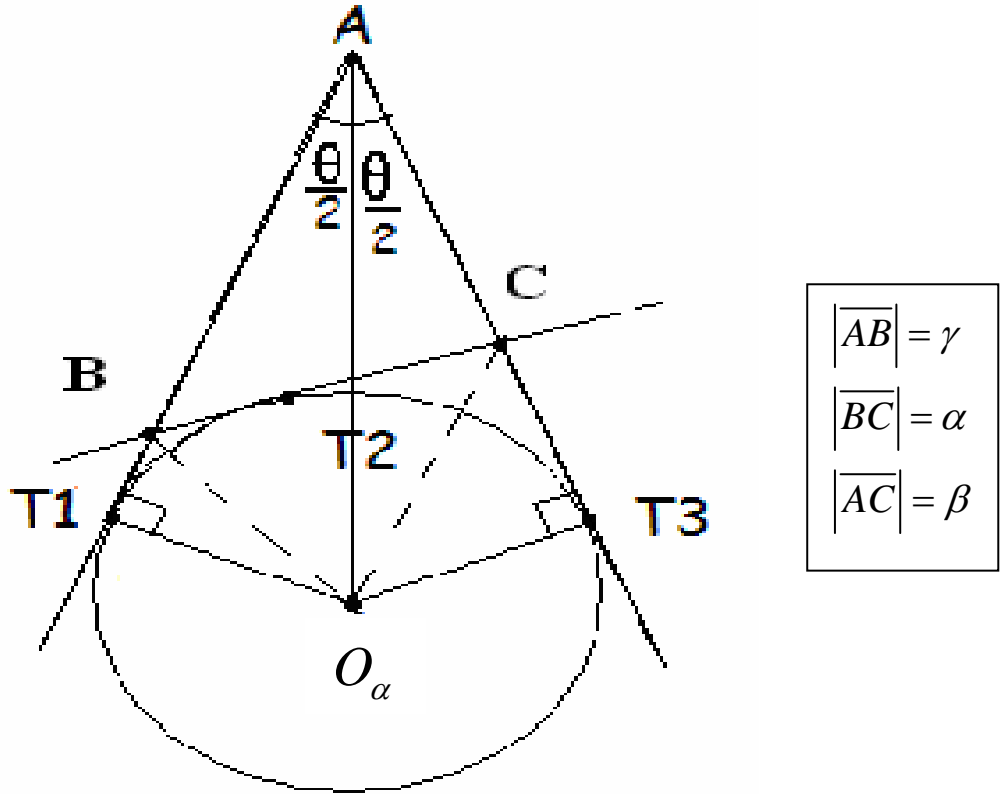

**Figure-1**

Throughout this paper, $s$ will stand for the triangle $ABC$'s half or semi perimeter: $2s = \alpha + \beta + \gamma \Leftrightarrow s = \dfrac{\alpha + \beta + \gamma}{2}$.

If we combine (1) with (2) and use $s = \dfrac{\alpha + \beta + \gamma}{2}$ we find that

$$\tan\left(\frac{\theta}{2}\right) = \sqrt{\frac{(s-\beta)(s-\gamma)}{s(s-\alpha)}}$$

(3)

# 3 Formulas for the external radii $\rho_\alpha, \rho_\beta, \rho_\gamma$

In the picture (Figure 1) we have the following:

$\left|\overline{AB}\right| = \gamma, \left|\overline{BC}\right| = \alpha, \left|\overline{AC}\right| = \beta, \left|\overline{O_\alpha T_1}\right| = \rho_\alpha = \left|\overline{O_\alpha T_3}\right|$ ; $T_1, T_2, T_3$ are the points of tangency. We put $x = \left|\overline{BT_1}\right| = \left|\overline{BT_2}\right|$ and $y = \left|\overline{T_2 C}\right| = \left|\overline{CT_3}\right|$.

Since $\left|\overline{AT_1}\right| = \left|\overline{AT_3}\right|$, it follows that $\gamma + x = \beta + y$; or equivalently

$$x - y = \beta - \gamma \qquad (4)$$

But also $\left|\overline{BT_2}\right| + \left|\overline{T_2 C}\right| = \left|\overline{BC}\right|$;

$$x + y = \alpha \qquad (5)$$

From (4) and (5) it follows that,

$$x = \frac{\beta - \gamma + \alpha}{2} \ , \ y = \frac{\alpha + \gamma - \beta}{2} \qquad (6)$$

and from the right angle $O_\alpha T_1 A$ we have

$$\rho_\alpha = \left|\overline{AT_1}\right| \cdot \tan\left(\frac{\theta}{2}\right) \qquad (7)$$

But $\left|\overline{AT_1}\right| = \gamma + x = \frac{\gamma + \beta + \alpha}{2} = s$; and also (7) gives

$$\rho_\alpha = s \tan\left(\frac{\theta}{2}\right) \qquad (8)$$

Combining (8) with (3) leads to the formula

$$\rho_\alpha = \sqrt{\frac{s(s - \beta)(s - \gamma)}{(s - \alpha)}} \qquad (9)$$

If $E$ is the area of a triangle $ABC$, the $E = \sqrt{s(s - \alpha)(s - \beta)(s - \gamma)}$ . This formula for the area $E$ in terms of the side lengths $\alpha, \beta, \gamma$ is known in the literature as Heron's formula.

If we use Heron's Formula and (9), we easily obtain the formula

$$\rho_\alpha = \frac{E}{s - \alpha},$$

and by cyclic interchange we also have

$$\rho_\alpha = \frac{E}{s - \beta} \ \text{ and } \ \rho_\alpha = \frac{E}{s - \gamma}$$

Note that when the triangle $AB\Gamma$ is a right-angled one with $\theta = 90$; then (8) easily implies $\rho_\alpha = s =$ semi perimeter.

We also have $\beta^2 + \gamma^2 = \alpha^2$ and area $E = \dfrac{\beta\gamma}{2}$. We invite the reader to do the necessary algebra in order to show that,

$$\rho_\beta = \frac{\beta - \gamma + \alpha}{2} \text{ and } \rho_\gamma = \frac{\alpha + \gamma - \beta}{2}. \tag{10}$$

The converse is also true: if $\rho_\alpha = s$, then $\theta = 90$. For a proof of the fact that $\rho_\alpha = s$ if, and only if, $\theta = 90$. The reader may refer to Amy Bell's paper (see [1]).

## 4 When $ABC$ is a Pythagorean triangle

A Pythagorean triangle is a right-angled Heron triangle. Basic material on Pythagorean triangles can be found in [5] and [12]. When $ABC$ is Pythagorean, then α, β, γ are positive integers such that $\alpha^2 = \beta^2 + \gamma^2$. This in turn implies that $\theta = 90$ and also,

$$\beta = \delta(2mn), \ \gamma = \delta(m^2 - n^2), \ \alpha = \delta(m^2 + n^2) \tag{11}$$

Or alternatively $\beta = \delta(m^2 - n^2)$, $\gamma = \delta(2mn)$; for some positive integers $m, n, \delta$ such that $m > n, (m, n) = 1$ (i.e., $m$ and $n$ are relatively prime) and $m + n \equiv 1 \pmod 2$ (One of $m$, $n$ is odd; the other even).

That the above parametric formulas describe the entire family of Pythagorean triples is well known. For historical information regarding Pythagorean triangles the reader may refer to [3]. Another classic book on the subject is [4] and for information on the subject from an elementary number theory view point see [5].

## 5 The External radii $\rho_\alpha, \rho_\beta, \rho_\gamma$ of a Pythagorean triangle

When $ABC$ is a Pythagorean triangle it may be shown from (10) and (11) that

$$\rho_\alpha = \delta m(n + m), \rho_\beta = \delta n(m + n), \rho_\gamma = \delta m(m - n) \tag{12}$$

# 6 The case of Heron isosceles triangles

In this section, we examine Heron's isosceles triangles. These are triangles with $\beta = \gamma \neq \alpha$ being integers and the area $E$ also being an integer. To be able to describe all such triangles, we make use of Proposition 1, proved below. To establish Proposition 1, we need the following result from number theory

*Result 1:* A positive integer $i$ is the $l^{th}$ ($l$ a positive integer) power of a positive rational number $r$ if and only if, $r$ is an integer. Equivalently, the $l^{th}$ root of a positive integer $i$ is either an integer or an irrational number. In particular, the square root $\sqrt{i}$ is either an integer or otherwise irrational.

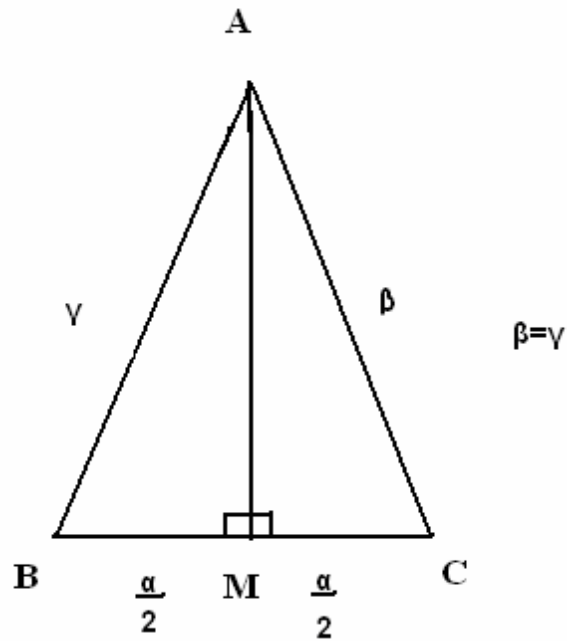

**Figure-2**

For more information and a proof of result 1, refer to [4] or [5].

**Proposition 1**: *Let $ABC$ be an isosceles triangle with side lengths*

$\gamma = \left| \overline{AB} \right| = \beta = \left| \overline{AC} \right| \neq \alpha = \left| \overline{BC} \right|$. *Suppose that* $\alpha$ *and* $\beta$ $(= \gamma)$ *are both (positive) integers and the area* $E$ *is also an integer. Also, let* $h$ *be the height from the vertex* $A$ *to the side* $\left| \overline{BC} \right|$. *Then* $\alpha$ *must be an even integer and* $h$ *an integer.*

**Proof:** From area $E = \dfrac{h.\alpha}{2} \Leftrightarrow h = \dfrac{2E}{\alpha}$ **,** it follows that $h$ is a rational number since both $E$ and $\alpha$ are integers. We put $h = \dfrac{M}{N}$ , where $M$ and $N$ are relatively prime positive integers. From either of the two congruent right triangles $ABM$ and $AMC\Gamma$ (Figure 2) it also follows that

$h^2 + \left( \dfrac{\alpha}{2} \right)^2 = \beta^2 \Leftrightarrow \left( \dfrac{M}{N} \right)^2 + \left( \dfrac{\alpha}{2} \right)^2 = \beta^2$ ; or equivalently,

$$4M^2 + \alpha^2 \cdot N^2 = 4N^2 \beta^2 \qquad (13)$$

Obviously, (13) shows that 4 must divide $\alpha^2 \cdot N^2 = (\alpha \cdot N)^2$, hence 2 must divide $\alpha N$. In other words, $\alpha \cdot N$ must be an integer, which means that at least one of $\alpha$ and $N$ must be even. If $N$ is even, $N = 2N_1$, for some positive integer $N_1$. Then by (13)

$$M^2 + \alpha^2 \cdot N_1^2 = 4N_1^2 \beta^2 \qquad (14)$$

Since $M$ and $N$ are relatively prime and $N$ is even, it follows that $M$ must be odd. But then, (14) easily implies that the integer $\alpha^2 \cdot N_1^2 = (\alpha \cdot N_1)^2$ must also be odd, which in turn implies that $\alpha N_1$ must be odd as well. Now, the square of any odd integer is congruent to 1 modulo 4 (in fact, 1 modulo 8, but we do not need that). Thus,

$$M^2 \equiv \left( \alpha N_1 \right)^2 \equiv 1 (\mathrm{mod}\, 4) \Rightarrow M^2 + \alpha^2 N_1^2 \equiv 2 (\mathrm{mod}\, 4),$$

and thus not a multiple of 4, therefore contradicting (14). Thus, $N$ cannot be even, and since $\alpha N$ is even, it follows that $\alpha$ is even. Since $\alpha$ is even $\alpha = 2\alpha_1$, for some positive integer $\alpha_1$. From $h^2 + \left( \dfrac{\alpha}{2} \right)^2 = \beta^2$ we obtain $h^2 = \beta^2 - \alpha_1^2$, a positive integer. Since $h$ is rational, it follows by Result 1, that in fact $h$ must be an integer.
The proof is complete.

The following proposition essentially describes the family of all heron isosceles triangles. In effect, proposition 2 shows that each Heron isosceles triangle can be obtained by glueing together two congruent Pythagorean triangles (so that they share a common leg).

**Proposition 2:** *Let $ABC$ be an isosceles triangle with $\gamma = \overline{|AB|} = \beta = \overline{|AC|} \neq \alpha = \overline{|BC|}$, the side lengths $\alpha, \beta, \gamma$ being integers and with the area $E$ also an integer. Also, let $h$ be the height from the vertex $A$ to the side $\overline{BC}$. Then, $\alpha = 2\delta(m^2 - n^2), h = \delta(2mn), \gamma = \beta = \delta(m^2 + n^2)$; Or alternatively $\alpha = 4\delta mn, h = \delta(m^2 - n^2), \gamma = \beta = \delta(m^2 + n^2)$, for some positive integers $\delta, m, n$ such that $m > n, (m,n) = 1$ ( i.e., $m$ and $n$ are relatively prime) and $m + n \equiv 1 \pmod 2$ (i.e., one of $m, n$ is odd; the other even). Conversely, if the side lengths $\alpha, \beta, \gamma$ of a triangle $ABC$ satisfy either set of formulas above for some positive integers $\delta, m, n$, such that $m > n$; then $h = \delta(2mn)$ or $h = \delta(m^2 - n^2)$ respectively, and $E$ is an integer.*

*Note: The conditions $(m,n) = 1$ and $m + n \equiv 1 \pmod 2$ are not really necessary for the converse statement.*

**Proof:** The converse part is straight forward computation and is left to the reader. Since $ABC$ is a Heron isosceles triangle with $\beta = \gamma \neq \alpha$. It follows from Proposition 1 that $\alpha$ is even and $h$ is an integer. Thus since $\frac{\alpha}{2}$ is an integer, both congruent right angled triangles $ABM$ and $AMC$ in Figure 2 are Pythagorean with hypotenuse length $\beta = \gamma$. Therefore, according to (11) we must have

$\alpha = 2\delta(m^2 - n^2), h = \delta(2mn), \gamma = \beta = \delta(m^2 + n^2)$ **;** or alternatively
$\alpha = 4\delta mn, h = \delta(m^2 - n^2), \gamma = \beta = \delta(m^2 + n^2)$ **,**

The proof is complete.

Now that we have a complete parametric description of the set of all Heron isosceles triangles, let us identify among such triangles, those for which the external radii $\rho_\alpha, \rho_\beta, \rho_\gamma$ are also integers. This is done in Theorem 1 below. To do so, however, we need another result from number theory as follows:

***Result 2:*** *let $a, b, c$ be integers. If $c$ is a divisor of the product $a.b$ and $c$ is relatively prime to $a$, then $c$ must be a divisor of $b$.*

Result 2 and a proof of it can be found, for example, in [4] or [5].

**Theorem 1:** *Let $S$ be the set of all Heron isosceles triangles; namely, the set of all triangles $ABC$ such that $\gamma = \overline{|AB|} = \beta = \overline{|AC|} \neq \alpha = \overline{|BC|}$, with the side lengths $\alpha, \beta, \gamma$ being integers and with the area $E$ also an integer; and with the added property that the three external radii $\rho_\alpha, \rho_\beta, \rho_\gamma$ are also integers. Then $S$ is the union of the following two families:*

*$F_1: \beta = \gamma = Kn(m^2 + n^2), \alpha = 2Kn(m^2 - n^2)$; in this family the radii are given by $\rho_\alpha = Km(m^2 - n^2), \rho_\beta = \rho_\gamma = 2Kmn^2$*

*$F_2: \beta = \gamma = L(m-n)(m^2 + n^2), \alpha = 4L(m-n)mn$;*                    and                    with  *$\rho_\alpha = 2Lmn(m+n), \rho_\beta = \rho_\gamma = L(m+n)(m-n)^2$*

*Where $K$ and $L$ are arbitrary positive integers; and $m, n$ can be any two relatively prime integers of different parities ( one even, the other odd), and with $m > n$.*

**Proof:** From the definition of the set $S$ and the parametric formulas that define the family $F_1$ and $F_2$, it is apparent that if a Heron triangle belongs in $F_1$ or $F_2$, it must belong in $S$. Below we prove the converse. That is, if a triangle is a member of $S$, it must be in $F_1$ or $F_2$. By Proposition 2, we know that if $ABC$ is a Heron isosceles triangle with $\beta = \gamma \neq \alpha$; we must have

$$\alpha = 2\delta(m^2 - n^2), h = \delta(2mn), \beta = \gamma = \delta(m^2 + n^2) \qquad (15a)$$

Or alternatively

$$\alpha = 4\delta mn., h = \delta(m^2 - n^2), \beta = \gamma = \delta(m^2 + n^2), \qquad (15b)$$

where $m, n$ are positive integers such that $(m,n) = 1$ , $m + n \equiv 1 \pmod{2}$, and $m > n$.

Also by (10) we know that

$$\rho_\alpha = \frac{E}{s - \alpha} \text{ and } \rho_\beta = \rho_\gamma = \frac{E}{s - \beta} \text{ since } (\beta = \gamma)$$

Since $E = \frac{\alpha}{2} \cdot h$ and $s = \frac{2\beta + \alpha}{2} = \beta + \frac{\alpha}{2}$, it follows that

$$\rho_\alpha = \frac{\alpha h}{2\beta - \alpha} \ \text{ and } \ \rho_\beta = \rho_\gamma = h \qquad (16)$$

If (15a) holds, (16) and a straight forward calculation yield

$$\rho_\alpha = \frac{\delta(m^2 - n^2)m}{n} \ \text{ and } \ \rho_\beta = \rho_\gamma = h = \delta(2mn) \qquad (17)$$

Since $m$ and $n$ are relatively prime, it follows that n must be relatively prime to the product $\{(m^2 - n^2)m\}$; this is left to the reader to verify. (Note that $(m,n) = 1$ alone implies $(n,(m^2 - n^2)m) = 1$; the condition $m + n \equiv 1 (\bmod\, 2)$ is not needed for this). On the other hand, as (16) clearly shows, $\rho_\alpha$ will be an integer if, and only if, $n$ is a divisor of the product $\delta(m^2 - n^2)m$. Therefore since $n$ is relatively prime to $(m^2 - n^2)m$, then by result 2, $n$ must be a divisor of $\delta$; $\delta = Kn$, for some positive integer $K$. Substituting for $\delta = Kn$ in (15a) and (17) produces

$$\alpha = 2Kn(m^2 - n^2), \beta = \gamma = kn(m^2 + n^2), \rho_\alpha = Km(m \ - n^2) \text{,and}$$
$$\rho_\beta = \rho_\gamma = h = 2Kmn^2$$

If on the other hand, possibility (15b) holds, then again combining (15b) with (16) gives,

$$\rho_\alpha = \frac{\delta(2mn)(m+n)}{m-n} \ \text{ and } \ \rho_\beta = \rho_\gamma = h = \delta(m^2 - n^2) \qquad (18)$$

Applying the same reasoning as before, we see that in order for $\rho_\alpha$ to be an integer, it is necessary and sufficient that $(m-n)$ be a divisor of $\delta$. This follows from Result 2 and the fact that the positive integer $(m-n)$ is relatively prime to the product $2mn(m+n)$. This last inference, in turn, follows from the conditions $(m,n) = 1$ and $m + n \equiv 1 (\bmod\, 2)$. We see that in this case, it is not enough that $m$ and $n$ are relatively prime; the fact that one of $m,n$ is even, while the other is odd, is also needed. We set $\delta = L(m-n)$, where $L$ is a positive integer. Combining (15b) with (18) yields

$$\beta = \gamma = L(m-n)(m^2+n^2), \alpha = 4L(m-n)mn \text{ , and}$$
$$\rho_\alpha = 2Lmn(m+n), \rho_\beta = \rho_\gamma = L(m+n)(m-n)^2 = h$$

The proof is complete.

# 7  Numerical examples

We present those Heron isosceles triangles with integral $\rho_\alpha, \rho_\beta, \rho_\gamma$, under the parameter constraints $1 \le n < m \le 6, L = K = 1, \quad (m,n)=1$, and $m+n \equiv 1 \pmod 2$,

Family $F_1$

|  | $\alpha$ | $\beta = \gamma$ | $\rho_\beta = \rho_\gamma$ | $\rho_\alpha$ |
|---|---|---|---|---|
| K=1, n=1, m=2 | 6 | 5 | 4 | 6 |
| K=1, n=1, m=4 | 30 | 17 | 8 | 60 |
| K=1, n=1, m=6 | 70 | 37 | 12 | 210 |
| K=1, n=2, m=3 | 20 | 26 | 24 | 15 |
| K=1, n=2, m=5 | 84 | 58 | 40 | 105 |
| K=1, n=3, m=4 | 42 | 75 | 72 | 28 |
| K=1, n=4, m=5 | 72 | 164 | 160 | 45 |
| K=1, n=5, m=6 | 110 | 305 | 300 | 66 |

Family $F_2$

|  | $\alpha$ | $\beta = \gamma$ | $\rho_\beta = \rho_\gamma$ | $\rho_\alpha$ |
|---|---|---|---|---|
| K=1, n=1, m=2 | 8 | 5 | 3 | 12 |
| K=1, n=1, m=4 | 48 | 51 | 45 | 40 |
| K=1, n=1, m=6 | 120 | 185 | 175 | 84 |
| K=1, n=2, m=3 | 24 | 13 | 5 | 60 |
| K=1, n=2, m=5 | 120 | 87 | 63 | 140 |

| K=1, n=3, m=4 | 48 | 25 | 7 | 168 |
|---|---|---|---|---|
| K=1, n=4, m=5 | 80 | 41 | 9 | 360 |
| K=1, n=5, m=6 | 120 | 61 | 11 | 660 |